# Block representations of generalized inverses of matrices

**Vera Miler Jerković**

*School of Electrical Engineering, University of Belgrade, 73 Bulevar kralja Aleksandra*
*e-mail: vera.miler@etf.rs*

**Branko Malešević**

*School of Electrical Engineering, University of Belgrade, 73 Bulevar kralja Aleksandra*
*e-mail: branko.malesevic@etf.rs*

**Abstract.** In this paper will be considered standard forms of generalized inverses for matrices in the shape of block representations $\{1, 2, 3, 4, 5, 5^k\}$-inverse. Especially will be considered Moore-Penrose inverse and the group inverse. Results from Rhode's technique are used and methods for calculating some inverse are shown on examples.

## 1. Introduction

The concept of generalized inverse for matrices was considered by E.H. Moore (1920.) and R. Penrose (1955.) [1], [2] and [3]. The generalized form of $\{1, 2, 3, 4, 5, 5^k\}$-inverse as well as some combination (the Moore-Penrose inverse, the group inverse, the Drazin inverse) play major role in solving problems in various areas of sciences, such as fuzzy mathematics, linear regression,... The combinations $\{1,3\}$ and $\{1,4\}$ of generalized inverses and also Moore-Penrose inverse have minimax properties and it can be applied in solving linear systems. Also, combination $\{1,2,3\}$ of generalized inverses can be used for finding least squared solution of linear system. The group inverse has many applications in singular differential equations; Markov chains iterate methods and so on. The generalized inverses have application in linear statistical modeling, especially in solving singularity of covariance matrix. Application of $\{2\}$-inverse is using in Newton methods for solving systems of nonlinear equations. The Drazin inverse is using for solving singular linear difference equations. The overview of applications can be found in textbooks [2] and [3]. Presentations of various forms of generalized inverse are given according to method of C.A. Rhode [4]. This method is also presented in [5] and [6]. The Rhode's method is applied in many fields. It can be used for finding solution of matrix equation $AXB = C$ where solution is described in the terms of the Rhode's general form of the $\{1\}$-inverse [6], [7], [8], [9], [10]; see also [11]. This method is also applicable for solving matrix equation $AXB = C$ using Penrose's general solution [12]. The aim of this paper is to describe generalized form of $\{1, 2, 3, 4, 5, 5^k\}$-inverse using Rhode's method. Consequently, we obtain forms of the Moore-Penrose inverse, the group inverse and the Drazin inverse. Using block representations of generalized inverse of matrix, we got detailed structure of generalized inverses. In this way, we can approach parts of generalized inverse and modify it to do better, which will our future plan. All relevant theorems and methods for pseudoinverses in this paper are presented by appropriate examples.

## 2. Block representations

For the matrix $A \in C^{m \times n}$ system of four Penrose's equations is consider:

$$AXA = A \tag{1}$$
$$XAX = X \tag{2}$$
$$(AX)^* = AX \tag{3}$$
$$(XA)^* = XA \tag{4}$$

where matrix $X \in C^{n \times m}$ is unknown. For square matrix, we add matrix equations:

$$AX = XA \tag{5}$$
$$A^k XA = A^k \tag{6}$$

where $ind(A) = k$ is index of the matrix $A$.

**Definition 1.** The index $k$ of matrix $A$ is smallest non-negative number such that the equality $rank(A^k) = rank(A^{k+1})$ is true.

Solutions of matrix equations (1), (2), (3), (4), (5) and (6) we defined as $\{1\}$, $\{2\}$, $\{3\}$, $\{4\}$, $\{5\}$ and $\{5^k\}$ generalized inverse of the matrix $A$.

Let $C_r^{m \times n}$ be a set of all matrices over set of complex numbers of order $m \times n$ with a rank $r$. For the matrix $A \in C_r^{m \times n}$ we can make expanded matrix $\begin{bmatrix} A & I_m \\ I_n & 0 \end{bmatrix}$ which can be transformed, by elementary transformations on the columns and rows, in equivalent matrix:

$$\begin{bmatrix} A & I_m \\ I_n & 0 \end{bmatrix} \sim \ldots \sim \begin{bmatrix} E_r & Q \\ P & 0 \end{bmatrix} \tag{7}$$

where $E_r \in C_r^{m \times n}$ is matrix with $r$ ones on first $r$ places of the main diagonal and zeros on the all other places. The matrices $Q \in C^{m \times m}$ and $P \in C^{n \times n}$ are regular and the following equality is true:

$$QAP = E_r = \begin{bmatrix} I_r & 0 \\ 0 & 0 \end{bmatrix}. \tag{8}$$

**Definition 2.** The generalized inverses of the matrix $A \in C_r^{m \times n}$, which satisfy some of matrix equations $(1) - (4)$, and also $(5) - (6)$ in the case of square matrix, can be defined as *block representations*:

$$X = P \cdot \begin{bmatrix} X_0 & X_1 \\ X_2 & X_3 \end{bmatrix} \cdot Q \tag{9}$$

where $X_0 \in C^{r \times r}, X_1 \in C^{r \times (m-r)}, X_2 \in C^{(n-r) \times r}, X_3 \in C^{(n-r) \times (m-r)}$ are appropriate submatrices.

In the next section are presented some essential theorems for block representations. Concretely, specifying theorems we describe block representations of $\{1, 2, 3, 4, 5, 5^k\}$ inverses. Likewise, using corollaries we describe combinations of these inverses, which are unique. The details of proofs of these theorems are shown in [2], [3], [4], [5] and [6].

**Theorem 1**. (**Generalized $\{1\}$ – inverse**) *For the matrix $A \in C_r^{m \times n}$, $r < min\{m, n\}$, let there be given regular matrices $Q \in C^{m \times m}$ and $P \in C^{n \times n}$ such that satisfy (8). The matrix X of the shape (9) satisfies matrix equation (1) if and only if*:

$$X = P \cdot \begin{bmatrix} I_r & X_1 \\ X_2 & X_3 \end{bmatrix} \cdot Q \tag{10}$$

where $X_1 \in C^{r \times (m-r)}, X_2 \in C^{(n-r) \times r}$, and $X_3 \in C^{(n-r) \times (m-r)}$ are arbitrary submatrices.

In the case when $m=r$, then the matrix $A$ is matrix of full rank by rows, hence submatrices $X_1$ and $X_3$ dissapear by dimension. Therefore, the matrix $X$ from (10) has the shape

$$X = P \cdot \begin{bmatrix} I_r \\ X_2 \end{bmatrix} \cdot Q \tag{11}$$

In the case when *n=r*, then the matrix *A* is matrix of full rank by columns, hence submatrices $X_2$ and $X_3$ dissapear by dimension. Therefore, the matrix *X* from (10) has the shape

$$X = P \cdot [I_r \quad X_1] \cdot Q \tag{12}$$

**Theorem 2. (Generalized {2} – inverse)** *For the matrix $A \in C_r^{m \times n}$, $r < \min\{m, n\}$, let there be given regular matrices $Q \in C^{m \times m}$ and $P \in C^{n \times n}$ such that satisfy (8). The matrix X of the shape (9) satisfies matrix equation (2) if and only if submatrices $X_0 \in C^{r \times r}$, $X_1 \in C^{r \times (m-r)}$, $X_2 \in C^{(n-r) \times r}$ and $X_3 \in C^{(n-r) \times (m-r)}$ satisfy matrix equations:*

$$X_0^2 = X_0, \; X_0 X_1 = X_1, \; X_2 X_0 = X_2, \; X_2 X_1 = X_3. \tag{13}$$

**Corollary 1. (Generalized {1, 2} – inverse)** *For the matrix $A \in C_r^{m \times n}$, $r < \min\{m, n\}$, let there be given regular matrices $Q \in C^{m \times m}$ and $P \in C^{n \times n}$ such that satisfy (8). The matrix X of the shape (9) satisfies matrix equations (1) and (2) if and only if:*

$$X = P \cdot \begin{bmatrix} I_r & X_1 \\ X_2 & X_2 X_1 \end{bmatrix} \cdot Q \tag{14}$$

*where $X_1 \in C^{r \times (m-r)}$ and $X_2 \in C^{(n-r) \times r}$ are arbitrary submatrices.*

In the case of *m=r* the matrix *A* is matrix of full rank by rows, hence submatrix $X_1$ dissapears by dimension. Therefore, the matrix *X* from (14) has the shape

$$X = P \cdot \begin{bmatrix} I_r \\ X_2 \end{bmatrix} \cdot Q \tag{15}$$

In the case of *n=r* the matrix *A* is matrix of full rank by columns, hence submatrix $X_2$ dissapears by dimension. Therefore, the matrix *X* from (14) has the shape

$$X = P \cdot [I_r \quad X_1] \cdot Q \tag{16}$$

To detect {3} and {4} inverses of matrix *A* in the shape of block, it is necessary to make square block matrices:

$$Q \cdot Q^* = \begin{bmatrix} S_1 & S_2 \\ S_3 & S_4 \end{bmatrix} \text{ and } P^* \cdot P = \begin{bmatrix} T_1 & T_2 \\ T_3 & T_4 \end{bmatrix} \tag{17}$$

with appropriate submatrices $S_1 \in C^{r \times r}, S_2 \in C^{r \times (m-r)}, S_3 \in C^{(m-r) \times r}, S_4 \in C^{(m-r) \times (m-r)}$ and $T_1 \in C^{r \times r}, T_2 \in C^{r \times (n-r)}, T_3 \in C^{(n-r) \times r}, T_4 \in C^{(n-r) \times (n-r)}$. Matrices $Q \cdot Q^* \in C^{m \times m}$ and $P^* \cdot P \in C^{n \times n}$ are Hermitian. According to [4] hold:

$$S_1^* = S_1, S_2^* = S_3, S_4^* = S_4, T_1^* = T_1, T_2^* = T_3, T_4^* = T_4 \tag{18}$$

and also square matrices $S_4$ and $T_4$ are invertible.

**Definition 3.** The matrix matrix $A \in C^{m \times m}$ is Hermitian if $A = A^*$.

**Theorem 3. (Generalized {3} – inverse)** *For the matrix $A \in C_r^{m \times n}$, $r < \min\{m, n\}$, let there be given regular matrices $Q \in C^{m \times m}$ and $P \in C^{n \times n}$ such that satisfy (8) and let square matrix $Q \cdot Q^*$ be in the shape (17). The matrix X of the shape (9) satisfies matrix equation (3) if and only if submatrices $X_0 \in C^{r \times r}$, and $X_1 \in C^{r \times (m-r)}$ satisfy:*

$$(S_1 - S_2 S_4^{-1} S_2^*) X_0^* = X_0 (S_1 - S_2 S_4^{-1} S_2^*) \text{ and } X_1 = -X_0 S_2 S_4^{-1} \tag{19}$$

*where $X_2 \in C^{(n-r) \times r}$ and $X_3 \in C^{(n-r) \times (m-r)}$ are arbitrary submatrices .*

**Corollary 2. (Generalized {1, 3} – inverse)** *For the matrix $A \in C_r^{m \times n}$, $r < \min\{m, n\}$, let there be given regular matrices $Q \in C^{m \times m}$ and $P \in C^{n \times n}$ such that satisfy (8) and let square matrix $Q \cdot Q^*$ be in the shape (17). Matrix X of the shape (9) satisfies matrix equations (1) and (3) if and only if:*

$$X = P \cdot \begin{bmatrix} I_r & -S_2 S_4^{-1} \\ X_2 & X_3 \end{bmatrix} \cdot Q \qquad (20)$$

*where $X_2 \in C^{(n-r) \times r}$ and $X_3 \in C^{(n-r) \times (m-r)}$ are arbitrary submatrices.*

In the case of $m=r$ the matrix $A$ is matrix of full rank by rows, hence submatrices $S_2$, $S_4$ and $X_3$ dissapear by dimension. Therefore, the matrix $X$ from (20) has the shape

$$X = P \cdot \begin{bmatrix} I_r \\ X_2 \end{bmatrix} \cdot Q \qquad (21)$$

In the case of $n=r$, then the matrix $A$ is matrix of full rank by columns, hence submatrices $X_2$ and $X_3$ dissapear by dimension. Therefore, the matrix $X$ from (20) has the shape

$$X = P \cdot [I_r \quad -S_2 S_4^{-1}] \cdot Q \qquad (22)$$

**Corollary 3. (Generalized {1, 2, 3} – inverse)** *For the matrix $A \in C_r^{m \times n}$, $r < \min\{m, n\}$, let there be given regular matrices $Q \in C^{m \times m}$ and $P \in C^{n \times n}$ such that satisfy (8) and let square matrix $Q \cdot Q^*$ be in the shape (17). The matrix X of the shape (9) satisfies matrix equations (1), (2) and (3) if and only if:*

$$X = P \cdot \begin{bmatrix} I_r & -S_2 S_4^{-1} \\ X_2 & X_2 \cdot (-S_2 S_4^{-1}) \end{bmatrix} \cdot Q \qquad (23)$$

*where $X_2 \in C^{(n-r) \times r}$ is arbitrary submatrix.*

In the case of $m=r$ the matrix $A$ is matrix of full rank by rows, hence submatrices $S_2$ and $S_4$ dissapear by dimension. Therefore, the matrix $X$ from (23) has the shape

$$X = P \cdot \begin{bmatrix} I_r \\ X_2 \end{bmatrix} \cdot Q \qquad (24)$$

In the case of $n=r$ the matrix $A$ is matrix of full rank by columns, hence submatrix $X_2$ dissapears by dimension. Therefore, the matrix $X$ from (23) has the shape

$$X = P \cdot [I_r \quad -S_2 S_4^{-1}] \cdot Q \qquad (25)$$

**Theorem 4. (Generalized {4} – inverse)** *For the matrix $A \in C_r^{m \times n}$, $r < \min\{m, n\}$, let there be given regular matrices $Q \in C^{m \times m}$ and $P \in C^{n \times n}$ such that satisfy (8) and let square matrix $P^* \cdot P$ be in the shape (17). Matrix X of the shape (9) satisfies matrix equation (4) if and only if submatrices $X_0 \in C^{r \times r}$ and $X_2 \in C^{(n-r) \times r}$ satisfy:*

$$X_0^*(T_1 - T_2 T_4^{-1} T_2^*) = (T_1 - T_2 T_4^{-1} T_2^*) X_0 \quad and \quad X_2 = -T_4^{-1} T_3 X_0 \qquad (26)$$

*where $X_0 \in C^{r \times r}$ and $X_2 \in C^{(n-r) \times r}$ are arbitrary submatrices.*

**Corollary 4. (Generalized {1, 4} – inverse)** *For the matrix $A \in C_r^{m \times n}$, $r < \min\{m, n\}$, let there be given regular matrices $Q \in C^{m \times m}$ and $P \in C^{n \times n}$ such that satisfy (8) and let square matrix $P^* \cdot P$ be in the shape (17). The matrix X of the shape (9) satisfies matrix equations (1) and (4) if and only if:*

$$X = P \cdot \begin{bmatrix} I_r & X_1 \\ -T_4^{-1} T_3 & X_3 \end{bmatrix} \cdot Q \qquad (27)$$

*where $X_1 \in C^{r \times (m-r)}$ and $X_3 \in C^{(n-r) \times (m-r)}$ are arbitrary submatrices.*

In the case of *m=r* the matrix *A* is matrix of full rank by rows, hence submatrices $X_1$ and $X_3$ dissapear by dimension. Therefore, the matrix *X* from (27) has the shape

$$X = P \cdot \begin{bmatrix} I_r \\ -T_4^{-1}T_3 \end{bmatrix} \cdot Q \qquad (28)$$

In the case of *n=r* the matrix *A* is matrix of full rank by columns, hence submatrices $T_3, T_4$ and $X_3$ dissapear by dimension. Therefore, the matrix *X* from (27) has the shape

$$X = P \cdot [I_r \quad X_1] \cdot Q \qquad (29)$$

**Corollary 5**. **(Generalized {1, 2, 4} – inverse)** *For the matrix $A \in C_r^{m \times n}$, $r < \min\{m, n\}$, let there be given regular $Q \in C^{m \times m}$ and $P \in C^{n \times n}$ such as to satisfy (8) and let square matrix $P^* \cdot P$ be in the shape (17). The matrix X of the shape (9) satisfies matrix equations (1), (2) and (4) if and only if:*

$$X = P \cdot \begin{bmatrix} I_r & X_1 \\ -T_4^{-1}T_3 & (-T_4^{-1}T_3) \cdot X_1 \end{bmatrix} \cdot Q \qquad (30)$$

*where $X_1 \in C^{r \times (m-r)}$ is arbitrary submatrix.*

In the case of *m=r* the matrix *A* is matrix of full rank by rows, hence submatrix $X_1$ dissapears by dimension. Therefore, the matrix *X* from (30) has the shape

$$X = P \cdot \begin{bmatrix} I_r \\ -T_4^{-1}T_3 \end{bmatrix} \cdot Q \qquad (31)$$

In the case of *n=r* the matrix *A* is matrix of full rank by columns, hence submatrices $T_3$ and $T_4$ dissapear by dimension. Therefore, the matrix *X* from (30) has the shape

$$X = P \cdot [I_r \quad X_1] \cdot Q \qquad (32)$$

**Corollary 6**. **(Generalized {1, 3, 4} – inverse)** *For the matrix $A \in C_r^{m \times n}$, $r < \min\{m, n\}$, let there be given regular matrices $Q \in C^{m \times m}$ and $P \in C^{n \times n}$ such that satisfy (8) and let square matrices $Q \cdot Q^*$ and $P^* \cdot P$ have the shape (17). The matrix X of the shape (9) satisfies matrix equations (1), (3) and (4) if and only if:*

$$X = P \cdot \begin{bmatrix} I_r & -S_2 S_4^{-1} \\ -T_4^{-1}T_3 & X_3 \end{bmatrix} \cdot Q \qquad (33)$$

*where $X_3 \in C^{(n-r) \times (m-r)}$ is arbitrary submatrix.*

In the case of *m=r* the matrix *A* is matrix of full rank by rows, hence submatrices $S_2, S_4$ and $X_3$ dissapear by dimension. Therefore, the matrix *X* from (33) has the shape

$$X = P \cdot \begin{bmatrix} I_r \\ -T_4^{-1}T_3 \end{bmatrix} \cdot Q \qquad (34)$$

In the case when *n=r*, then the matrix *A* is matrix of full rank by columns, hence submatrices $T_3, T_4$ and $X_3$ dissapear by dimension. Therefore, the matrix *X* from (33) has the shape

$$X = P \cdot [I_r \quad -S_2 S_4^{-1}] \cdot Q \qquad (35)$$

For the matrix $A \in C_r^{m \times n}$ system of matrix equations (1), (2), (3) and (4) has unique solution that is denoted by $A^\dagger$ and known as *Moore-Penrose inverse of matrix*. From above theorems we presente the block representation of Moore-Penrose inverse:

**Theorem 5. (Moore-Penrose inverse)** *For the matrix $A \in C_r^{m \times n}$, $r < \min\{m,n\}$, let there be given regular matrices $Q \in C^{m \times m}$ and $P \in C^{n \times n}$ such that satisfy (8) and let square matrices $Q \cdot Q^*$ and $P^* \cdot P$ have the shape (17). Unique solution of matrix equations (1), (2), (3) and (4) is given with*

$$A^\dagger = P \cdot \begin{bmatrix} I_r & -S_2 S_4^{-1} \\ -T_4^{-1} T_3 & T_4^{-1} T_3 S_2 S_4^{-1} \end{bmatrix} \cdot Q \qquad (36)$$

In the case of $m=r$ the matrix $A$ is matrix of full rank by rows, hence submatrices $S_2$ and $S_4$ dissapear by dimension. Therefore, the matrix $A^\dagger$ from (36) has the shape

$$A^\dagger = P \cdot \begin{bmatrix} I_r \\ -T_4^{-1} T_3 \end{bmatrix} \cdot Q \qquad (37)$$

In the case of $n=r$ the matrix $A$ is matrix of full rank by columns, hence submatrices $T_3$ and $T_4$ dissapear by dimension. Therefore, the matrix $A^\dagger$ from (36) has the shape

$$A^\dagger = P \cdot [I_r \quad -S_2 S_4^{-1}] \cdot Q \qquad (38)$$

In the case that the matrix $A$ is regular and square, $m=n=r$, the Moore-Penrose inverse $A^\dagger$ of the matrix $A$ is equal with $A^{-1}$. In that case, the matrix $A^\dagger$ from (36) has the shape

$$A^\dagger = A^{-1} = P \cdot Q \qquad (39)$$

In next example, determination of the Moore-penrose inverse of matrix in a shape of block is presented.

**Example 1**. *For the matrix* $A = \begin{bmatrix} 1 & 2 & 3 \\ 4 & 5 & 6 \\ 7 & 8 & 9 \end{bmatrix}$ *calculate the Moore-Penrose inverse.*

*Solution.* Using elementary transformations by columns and rows of expanded matrix $\begin{bmatrix} A & I_3 \\ I_3 & 0 \end{bmatrix}$ we can determine regular matrices $P$ and $Q$:

$$P = \begin{bmatrix} 1 & 0 & 1 \\ 0 & 1 & -2 \\ 0 & 0 & 1 \end{bmatrix} \text{ and } Q = \begin{bmatrix} -\frac{5}{3} & \frac{2}{3} & 0 \\ \frac{4}{3} & -\frac{1}{3} & 0 \\ 1 & -2 & 1 \end{bmatrix}$$

such that $QAP = E_2$. By making product of matrices

$$Q \cdot Q^* = \begin{bmatrix} \frac{29}{9} & -\frac{22}{9} & -3 \\ -\frac{22}{9} & \frac{17}{9} & 2 \\ -3 & 2 & 6 \end{bmatrix} \text{ and } P^* \cdot P = \begin{bmatrix} 1 & 0 & 1 \\ 0 & 1 & -2 \\ 1 & -2 & 6 \end{bmatrix}$$

we obtain submatrices $S_2 = \begin{bmatrix} -3 \\ 2 \end{bmatrix}$, $S_4 = [6]$, $T_3 = [1 \quad -2]$, $T_4 = [6]$ based on which we can determine submatrices $-S_2 S_4^{-1} = \begin{bmatrix} \frac{1}{2} \\ -\frac{1}{3} \end{bmatrix}$, $-T_4^{-1} T_3 = \begin{bmatrix} -\frac{1}{6} & \frac{1}{3} \end{bmatrix}$, $(T_4^{-1} T_3)(S_2 S_4^{-1}) = \begin{bmatrix} -\frac{7}{36} \end{bmatrix}$. Therefore, the Moore-Penrose inverse of matrix $A$, according to theorem 5, is given by

$$A^\dagger = P \cdot \begin{bmatrix} 1 & 0 & \frac{1}{2} \\ 0 & 1 & -\frac{1}{3} \\ -\frac{1}{6} & \frac{1}{3} & -\frac{7}{36} \end{bmatrix} \cdot Q = \begin{bmatrix} -\frac{23}{36} & -\frac{1}{6} & \frac{11}{36} \\ -\frac{1}{18} & 0 & \frac{1}{18} \\ \frac{19}{36} & \frac{1}{6} & -\frac{7}{36} \end{bmatrix}$$

∆

**Block representations of generalized inverses of squared matrices**

In this part of paper shall be considered block representations of generalized inverses of the square matrix $A \in C_r^{n \times n}$. Let the matrix $A$ be square matrix with index $ind(A) = k$, and let given minimal polynomial be $\mu(x)$. Then:

$$\mu(x) = x^m + c_{m-1}x^{m-1} + \cdots + c_k x^k, \quad c_k \neq 0. \tag{40}$$

Using to minimal polynomial we can form $q$–polynomial:

$$q(x) = -\frac{1}{c_k} \cdot (x^{m-k-1} + c_{m-1} x^{m-k-2} + \cdots + c_{k+1}). \tag{41}$$

Relation between minimal polynomial (40) and $q$ – polynomial (41) can be defined by:

$$\mu(x) = c_k x^k (1 - xq(x)) \tag{42}$$

**Definition 4.** Let is matrix $A \in C_r^{n \times n}$ with index $ind(A) \leq 1$. Then system of matrix equations (1), (2) and (5) has unique solution denoted by $A^\#$ and known as *the group inverse*.

If the matrix $A$ has $ind(A) = 0$, then it is regular and hold $A^\# = A^{-1}$.

**Theorem 6. (Group inverse)** *For the square matrix $A \in C_r^{n \times n}$ with index $ind(A) = 1$ and adequate q–polynomial $q(x)$, the matrix $q(A)$ represents one {1}-inverse of matrix A. Unique solution of system of matrix equations* (1), (2) *and* (5) *given by q–polynomial is*

$$A^\# = A(q(A))^2. \tag{43}$$

Next theorem is shown and proved in [13]. According to block representation of the group inverse of square matrix can be obtained without determining $q$–polynomial.

**Theorem 7. (Group inverse)** *For square matrix $A \in C_r^{n \times n}$, $r < \min\{m, n\}$, with index $ind(A) = 1$ let there be given regular matrices $Q, P \in C^{n \times n}$ such that satisfy* (8). *Let block decomposition*

$$Q \cdot P = \begin{bmatrix} V_1 & V_2 \\ V_3 & V_4 \end{bmatrix} \tag{44}$$

*is true under assumption that the submatrix $V_4 \in C^{(n-r) \times (n-r)}$ is regular. Unique solution of system of matrix equations* (1), (2) *and* (5) *is given by block representation*

$$A^\# = P \cdot \begin{bmatrix} I_r & -V_2 V_4^{-1} \\ -V_4^{-1} V_3 & V_4^{-1} V_3 V_2 V_4^{-1} \end{bmatrix} \cdot Q \tag{45}$$

The theorems 6 and 7 are illustrated by next example.

**Example 2**. *For the matrix* $A = \begin{bmatrix} 1 & 2 & 3 \\ 4 & 5 & 6 \\ 7 & 8 & 9 \end{bmatrix}$ *calculate the group inverse.*

*Solution.* Minimal polynomial and $q$–polynomial of the matrix $A$ are $\mu(x) = x^3 - 15x^2 - 18x$ and $q(x) = \frac{x}{18} - \frac{5}{6}$ respecivly. The index of matrix $A$ is $ind(A) = 1$. The group inverse of the matrix $A$, according to theorem 6, is

$$A^\# = A(q(A))^2 = \begin{bmatrix} -\frac{23}{36} & -\frac{1}{6} & \frac{11}{36} \\ -\frac{1}{18} & 0 & \frac{1}{18} \\ \frac{19}{36} & \frac{1}{6} & -\frac{7}{36} \end{bmatrix}$$

Above result can be reached using theorem 7. By use of the elementary transformation on columns and rows of expanded matrix $\begin{bmatrix} A & I_3 \\ I_3 & 0 \end{bmatrix}$ we can determine regular matrices $P$ and $Q$:

$$P = \begin{bmatrix} 1 & 0 & 1 \\ 0 & 1 & -2 \\ 0 & 0 & 1 \end{bmatrix}, \quad Q = \begin{bmatrix} -\frac{5}{3} & \frac{2}{3} & 0 \\ \frac{4}{3} & -\frac{1}{3} & 0 \\ 1 & -2 & 1 \end{bmatrix}$$

such that $QAP = E_2$. From product of matrices

$$Q \cdot P = \begin{bmatrix} -\frac{5}{3} & \frac{2}{3} & -3 \\ \frac{4}{3} & -\frac{1}{3} & 2 \\ 1 & -2 & 6 \end{bmatrix}$$

we can determine submatrices: $V_1 = \begin{bmatrix} -\frac{5}{3} & \frac{2}{3} \\ \frac{4}{3} & -\frac{1}{3} \end{bmatrix}$, $V_2 = \begin{bmatrix} -3 \\ 2 \end{bmatrix}$, $V_3 = [1 \ -2]$, $V_4 = [6]$. Next, according to $|V_4| = 6 \neq 0$, we get

$$A^{\#} = P \cdot \begin{bmatrix} I_r & -V_2 V_4^{-1} \\ -V_4^{-1} V_3 & V_4^{-1} V_3 V_2 V_4^{-1} \end{bmatrix} \cdot Q = P \cdot \begin{bmatrix} 1 & 0 & \frac{1}{2} \\ 0 & 1 & -\frac{1}{3} \\ -\frac{1}{6} & \frac{1}{3} & -\frac{7}{36} \end{bmatrix} \cdot Q = \begin{bmatrix} -\frac{23}{36} & -\frac{1}{6} & \frac{11}{36} \\ -\frac{1}{18} & 0 & \frac{1}{18} \\ \frac{19}{36} & \frac{1}{6} & -\frac{7}{36} \end{bmatrix}$$

Given solutions according to theorem 6 and theorem 7 are the same. △

**Definition 5.** The square matrix $A \in C_r^{n \times n}$ with index $ind(A) = k$, system of matrix equations (2), (5) and (6) has unique solution denoted by $A^D$ and known as *the Drazin inverse of matrix*.

**Theorem 8. (Drazin inverse)** *For the square matrix $A \in C_r^{n \times n}$ with index $ind(A) = k$ and adequate q–polynomial, unique solution of system of matrix equations (2), (5) and (6) is given with*
$$A^D = A^k (q(A))^{k+1} \tag{46}$$

If $k = ind(A) \leq 1$, then $A^D = A^{\#}$.

**Theorem 9.** [3] *For the square matrix $A \in C_r^{n \times n}$ holds $AA^D A = A$ if and only if $ind(A) \leq 1$.*

**Definition 6.** The square matrix $A \in C_r^{n \times n}$ for which holds true that $A^{\dagger} = A^D$ is called *EP–matrix*.

According to theorem 8, for *EP*–matrices holds $A^{\dagger} = A^D = A^{\#}$. At the end of this section, next theorem determines one characterisation of *EP*–matrices:

**Theorem 10.** [13] *Square matrix $A \in C_r^{n \times n}$ is EP–matrix if and only if one of equivalent condition satisfied*:
1. $N(A) = N(A^*)$.
2. $R(A) = R(A^*)$.
3. $\mathbf{C}^n = N(A) \oplus R(A)$.

4. There exist regular matrices $P \in C^{n \times n}$ and $A_r \in C^{r \times r}$ such that:

$$A = P \cdot \begin{bmatrix} A_r & 0 \\ 0 & 0 \end{bmatrix} \cdot P^{-1} \qquad (47)$$

**Corollary 6.** Every square singular symmetrical matrix $A \in C_r^{n \times n}$ is EP–matrix.

Based on previous two theorems here we give one example of computing Moore-Penrose inverse of EP– matrix.

**Example 3.** *Calculate the Moore-Penrose inverse of matrix*:

$$A = \begin{bmatrix} 1 & 1 & 1 & 0 & 0 \\ 1 & 2 & 0 & 1 & 1 \\ 1 & 0 & 2 & -1 & -1 \\ 0 & 1 & -1 & 1 & 1 \\ 0 & 1 & -1 & 1 & 1 \end{bmatrix}$$

*Solution.* Let us note that the matrix $A$ is symmetrical and singular ($|A| = 0$). It is meaning that matrix $A$ is one EP–matrix with index $ind(A) = 1$. So, $A^\dagger = A^D = A^\#$, and we can applied theorem 7. By applying elementary transformations on columns and rows of expanded matrix $\begin{bmatrix} A & I_5 \\ I_5 & 0 \end{bmatrix}$ we can obtain regular matrices $P$ and $Q$:

$$P = \begin{bmatrix} 1 & 0 & -2 & 1 & 1 \\ 0 & 1 & 1 & -1 & -1 \\ 0 & 0 & 1 & 0 & 0 \\ 0 & 0 & 0 & 1 & 0 \\ 0 & 0 & 0 & 0 & 1 \end{bmatrix} \text{ and } Q = \begin{bmatrix} 2 & -1 & 0 & 0 & 0 \\ -1 & 1 & 0 & 0 & 0 \\ -2 & 1 & 1 & 0 & 0 \\ 1 & -1 & 0 & 1 & 0 \\ 1 & -1 & 0 & 0 & 1 \end{bmatrix}$$

such that $QAP = E_2$. From the product of matrices

$$Q \cdot P = \begin{bmatrix} 2 & -1 & -5 & 3 & 3 \\ -1 & 1 & 3 & -2 & -2 \\ -2 & 1 & 6 & -3 & -3 \\ 1 & -1 & -3 & 3 & 2 \\ 1 & -1 & -3 & 2 & 3 \end{bmatrix}$$

we can determine the submatrices $V_1 = \begin{bmatrix} 2 & -1 \\ -1 & 1 \end{bmatrix}$, $V_2 = \begin{bmatrix} -5 & 3 & 3 \\ 3 & -2 & -2 \end{bmatrix}$, $V_3 = \begin{bmatrix} -2 & 1 \\ 1 & -1 \\ 1 & -1 \end{bmatrix}$ and

$V_4 = \begin{bmatrix} 6 & -3 & -3 \\ -3 & 3 & 2 \\ -3 & 2 & 3 \end{bmatrix}$. Then, according to $|V_4| = 12 \neq 0$, we get

$$A^\dagger = A^\# = P \cdot \begin{bmatrix} I_r & -V_2 V_4^{-1} \\ -V_4^{-1} V_3 & V_4^{-1} V_3 V_2 V_4^{-1} \end{bmatrix} \cdot Q = P \cdot \begin{bmatrix} 1 & 0 & \frac{7}{12} & -\frac{1}{4} & -\frac{1}{4} \\ 0 & 1 & \frac{1}{4} & -\frac{1}{4} & \frac{1}{4} \\ \frac{1}{3} & \frac{1}{12} & \frac{25}{144} & -\frac{1}{16} & -\frac{1}{16} \\ 0 & \frac{1}{4} & -\frac{1}{16} & \frac{1}{16} & \frac{1}{16} \\ 0 & \frac{1}{4} & -\frac{1}{16} & \frac{1}{16} & \frac{1}{16} \end{bmatrix} \cdot Q,$$

and

$$A^\dagger = \begin{bmatrix} \frac{1}{9} & \frac{1}{9} & \frac{1}{9} & 0 & 0 \\ \frac{1}{9} & \frac{25}{144} & \frac{7}{144} & \frac{1}{16} & \frac{1}{16} \\ \frac{1}{9} & \frac{7}{144} & \frac{25}{144} & -\frac{1}{16} & -\frac{1}{16} \\ 0 & \frac{1}{16} & -\frac{1}{16} & \frac{1}{16} & \frac{1}{16} \\ 0 & \frac{1}{16} & -\frac{1}{16} & \frac{1}{16} & \frac{1}{16} \end{bmatrix}$$

△

## 3. Conclusion

In this paper were presented theorems for obtaining Moore-Penrose inverse of matrices based on block representations. In addition, one part of this paper was dedicate to EP matrices. Clases {1,3} and {1,4} of generalized inverses and Moore-Penrose inverse have some minimax properties which can be applied on linear systems. Described procedure for detecting block representations of generalized inverses of matrices can be used in image reconstruction [14]. In the future, this will be subject of the further research.

**Acknowledgment.** *This work was supported in part by the Serbian Ministry of Education, Science and Technological Development under Grant No. 175016 and Grants No. 174032 and 44006.*

## References

[1] **R. Penrose.** A generalized inverses for matrices. *Mathematical Proceedings of the Cambridge Philosophical Society* 1955, 51, 406-413

[2] **A. Ben-Israel ,T.N.E. Greville.** Generalized Inverses, Theory and Applications. *Springer, New York*, 2003.

[3] **S.L. Cambell, C.D. Meyer**. Inverses of Linear Transformations. *Siam, Philadelphia*, 2009.

[4] **C.A. Rhode**. Contribution to the theory, computation and application of generalized inverses (PhD dissertation). University of North Carolina at Releigh, 1964.

[5] **V. Perić.** Generalizirana reciproka matrice, *Stručno-metodički časopis Matematika*, Zagreb, 1982, 11 (1), 40-57.

[6] **B. Malešević.** Grupna funkcionalna jednačina (magistarski rad). *Matematički fakultet, Univerzitet u Beogradu*, 1998.

[7] **I. Jovović, B. Malešević.** A note on Solutions of the Matrix Equations AXB=C. *Scientific publications of the state university of Novi Pazar, Ser.A: Appl. Math. Inform. And Mech.* 2014, 6 (1),  45-55

[8] **B. Malešević, B. Radičić.** Reproductive and non-reproductive solutions of matrix equation AXB = C. *Proceedings of the second symposium on Matematics and Applications Faculty of Mathematics*, 2011, 157-163

[9] **B. Malešević, B. Radičić.** Non-reproductive and reproductive solutions of some matrix equations. *Proceedings of the Conference on Mathematical and Informational Technologies MIT,* 2011, 246-251

[10] **B. Malešević, B. Radičić.** Some considerations of matrix equations using the concept of reproductivity. *Kragujevac Journal of Mathematics*, 2012, 36 (1), 151-161

[11] **B. Malešević, I. Jovović, M. Makragić, B. Radičić.** A Note on Solutions of Linear Systems. *ISRN Algebra*, vol. 2013 Article ID 142124, 6 pages,2013.doi:10.1155/2013/142124

[12] **B. Radičić, B. Malešević.** Some considerations in Relation to the Matrix Equation AXB=C. *The Mediterranean Journal of Mathematics* 2014, 11, 841-856

[13] **N. Matzakos, D. Pappas**. EP matrices: Computation of the Moore-Penrose inverse via factorizations. *Journal of Applied Mathematics and Computing*, 2010, 34, 113-127

[14] **B. Malešević, B. Banjac, M. Makragić, R. Obradović**. Application of polynomial texture mapping in process of digitalization of cultural heritage, *arXiv:1312.6935* [*cs.GR*]